\documentclass{amsart}
\usepackage{amsmath,amssymb,amsthm}
\usepackage{tikz}
\usepackage[hidelinks]{hyperref}

\DeclareMathOperator{\reg}{reg}
\DeclareMathOperator{\lcm}{lcm}

\newtheorem{Theorem}{Theorem}[section]
\newtheorem{Lemma}[Theorem]{Lemma}
\newtheorem{Corollary}[Theorem]{Corollary}
\newtheorem{Proposition}[Theorem]{Proposition}
\newtheorem{IntroTheorem}{Theorem}

\theoremstyle{definition}
\newtheorem{Definition}[Theorem]{Definition}
\newtheorem{Example}[Theorem]{Example}

\theoremstyle{remark}
\newtheorem{Remark}[Theorem]{Remark}

\begin{document}

\title[Squarefree powers of simplicial forests]{Squarefree powers of simplicial forests}

\author[M. Amalore Nambi]{Marie Amalore Nambi}
\address{Sabanci University, Faculty of Engineering and Natural Sciences, Orta Mahalle, Tuzla 34956, Istanbul, Turkey}
\email{amalore.pushparaj@sabanciuniv.edu}

\author[F. Navarra]{Francesco Navarra}
\address{Sabanci University, Faculty of Engineering and Natural Sciences, Orta Mahalle, Tuzla 34956, Istanbul, Turkey}
\email{francesco.navarra@sabanciuniv.edu}

\author[A. A. Qureshi]{Ayesha Asloob Qureshi}
\address{Sabanci University, Faculty of Engineering and Natural Sciences, Orta Mahalle, Tuzla 34956, Istanbul, Turkey}
\email{aqureshi@sabanciuniv.edu, ayesha.asloob@sabanciuniv.edu}

\author[D. Veer]{Dharm Veer}
\address{Department of Mathematics, Purdue University, 150 N. University Street, West Lafayette, IN 47907-2067, USA}
\email{dveer@purdue.edu}


\begin{abstract}

The highest nonzero squarefree power of an edge ideal
always has linear quotients, whereas this property
can fail even for the facet ideal of a pure
two-dimensional simplicial tree. This raises the question of which combinatorial conditions ensure linear quotients beyond the graph case. We characterize the squarefree powers of facet ideals of pure simplicial forests that have linear quotients, showing that this property is equivalent to having a linear resolution or being linearly related. We also investigate the corresponding question for undirected path ideals of forests. Here the highest nonzero squarefree power always has linear quotients for \(3\)-path ideals, whereas the analogous assertion fails for \(4\)-path ideals of trees.
\end{abstract}
\maketitle

\section{Introduction}
Let $S=K[x_1,\ldots,x_n]$ be a polynomial ring over a field $K$, and
let $I\subseteq S$ be a squarefree monomial ideal. The
{\em $k$-th squarefree power} $I^{[k]}$ is the ideal generated by the
squarefree monomials in $I^k$. When $I=I(\Delta)$ is the facet ideal
of a simplicial complex $\Delta$, these powers are closely related
to matchings: $I(\Delta)^{[k]}$ is generated by the monomials whose
supports are unions of $k$ pairwise disjoint facets. In particular,
$I(\Delta)^{[k]}\ne0$ if and only if $1\leq k\leq\nu(\Delta)$,
where $\nu(\Delta)$ is the matching number of $\Delta$. 

For the edge ideal $I(G)$ of a finite simple graph $G$, the highest
nonzero squarefree power always has linear quotients. This reflects a classical fact from matching theory:
the vertex sets covered by maximum matchings form the bases of a
matroid~\cite[Sections~1 and~6]{EF}. Hence $I(G)^{[\nu(G)]}$ is
matroidal and has linear quotients~\cite[Theorem~12.6.2]{HHmonomialideals};
a direct proof was given in~\cite[Theorem~5.1]{BHZ18powersedgeideals}.

This raises the question of whether the same property holds
for facet ideals of higher-dimensional simplicial complexes.
Since the answer is negative in general, it is natural to
seek well-structured classes of simplicial complexes for
which this property might still hold. A natural candidate is the class of simplicial forests, introduced by Faridi~\cite{Faridi} as a generalization of forests to simplicial complexes. Squarefree powers of facet ideals of simplicial forests were studied systematically in~\cite{KNQ24squarefree,DRS}. However, even for simplicial forests, the highest nonzero squarefree power need not have a linear resolution. Indeed,
\cite[Example~2.3]{KNQ24squarefree} provides such a
counterexample for a pure two-dimensional simplicial tree.
Moreover, for a pure simplicial forest of positive dimension,
only the highest and next-to-highest nonzero squarefree powers
can have linear resolutions
\cite[Corollary~3.4]{KNQ24squarefree}.
The problem is therefore to characterize precisely when
these remaining powers have linear resolutions or linear
quotients.

For edge ideals, matching-theoretic conditions have been
used to characterize linear resolutions of squarefree
powers of forests~\cite{EHforests} and linear first
syzygies of squarefree powers of arbitrary graphs
\cite[Theorem~3.1]{NQT26}.
Motivated by these results, we obtain a combinatorial
characterization for pure simplicial forests and show
that having linear first syzygies is equivalent to
having linear quotients. Recall that an ideal generated
in one degree is {\em linearly related} if its first
syzygy module is generated by linear relations. Our first main result characterizes the linearly related squarefree
powers and proves that each of them has linear quotients and a linear
resolution.

\begin{IntroTheorem}[Theorem~\ref{thm:main}]\label{intro:main}
Let $\Delta$ be a pure simplicial forest of dimension $d\geq1$, and
let $1\leq p\leq\nu(\Delta)$. The following conditions are equivalent:
\begin{enumerate}
\item[\textnormal{(a)}] $I(\Delta)^{[p]}$ has linear quotients;
\item[\textnormal{(b)}] $I(\Delta)^{[p]}$ has a linear resolution;
\item[\textnormal{(c)}] $I(\Delta)^{[p]}$ is linearly related;
\item[\textnormal{(d)}] $\Delta$ satisfies the conditions
\textnormal{(C1)$_p$} and \textnormal{(C2)$_p$} of
Definition~\ref{def:conditions}.
\end{enumerate}
\end{IntroTheorem}

The conditions concern proper chains of facets and exchanges among
matchings, so the characterization is independent of $K$.
For $p=\nu(\Delta)$, condition (C2)$_p$ is vacuous
(Corollary~\ref{cor:top}).


Directed path ideals of rooted trees provide an application of
Theorem~\ref{intro:main};  for background see~\cite{ConcaDeNegri,JV}.
Let $\Gamma$ be a rooted tree with edges directed away from the
root, and let $\Gamma_t$ be the simplicial complex of its directed
$t$-vertex paths. Since $\Gamma_t$ is a pure simplicial tree
whenever nonempty~\cite[Theorem~2.7]{JV}, our theorem applies.
Corollary~\ref{cor:rooted-trees} translates its conditions into
branching and connectivity properties of $\Gamma$, extending the
highest-power result for broom graphs
in~\cite[Theorem~2.7]{KNQ24squarefree}.

This naturally raises the question of what happens when all
paths of a forest are considered, rather than only directed
paths of a rooted tree. For a graph $G$, let $G^{(t)}$ be the
hypergraph of its $t$-vertex paths and set
$I_t(G)=I(G^{(t)})$. Since $I_2(G)=I(G)$, the first new case
is $t=3$. Although the associated path complex of a forest
need not be a simplicial forest, we prove that the highest
nonzero squarefree power of $I_3(G)$ has linear quotients.

\begin{IntroTheorem}[Theorem~\ref{thm:3path}]\label{intro:3path}
Let $G$ be a forest, and let $\nu=\nu(G^{(3)})\geq1$.
Then $I_3(G)^{[\nu]}$ has linear quotients. In particular, it has a
$3\nu$-linear resolution.
\end{IntroTheorem}


Moreover, we provide examples showing that the linear resolution conclusion fails for $4$-path ideals, even for trees (Example~\ref{ex:four-paths}), and for $3$-path ideals of arbitrary graphs (Example~\ref{ex.arbG}).

The paper is organized as follows. Section~\ref{sec:prelim} contains
the preliminary results on simplicial forests and matchings.
Section~\ref{sec:linear-related} establishes the characterization of
linearly related squarefree powers, and Section~\ref{sec:LQ} completes
the proof of Theorem~\ref{intro:main} by constructing linear quotients. 
Section~\ref{sec:path-ideals} contains the application to directed
path ideals of rooted trees and the proof of Theorem~\ref{intro:3path}.

\section{Simplicial forests, squarefree powers, and matchings}\label{sec:prelim}


Let $\Delta$ be a simplicial complex on a finite vertex set $V$.
We write $\mathcal F(\Delta)$ for its set of facets and
$\Delta=\langle F_1,\ldots,F_m\rangle$ when
$\mathcal F(\Delta)=\{F_1,\ldots,F_m\}$.
The complex is {\em pure} of dimension $d$ if all its facets have
$d+1$ vertices.
A {\em subcollection} of $\Delta$ is a simplicial complex whose facets
are facets of $\Delta$. For $W\subseteq V$, the {\em induced subcollection}
$\Delta_W$ is the complex generated by the facets of $\Delta$ contained
in $W$; it is not, in general, the restriction of $\Delta$ to $W$.
We assign matching number~$0$ to the empty subcollection.
Let $S=K[x_v:v\in V]$ be a polynomial ring over a field $K$, and set
$x_W=\prod_{v\in W}x_v$, with $x_{\varnothing}=1$.
The {\em facet ideal} of $\Delta$ is
$I(\Delta)=(x_F:F\in\mathcal F(\Delta))$.

A facet $F$ of $\Delta$ is a {\em leaf} if it is the only facet, or if
there is a facet $G\ne F$, called a {\em joint} of $F$, such that
$F\cap H\subseteq F\cap G$ for every facet $H\ne F$.
The complex is a {\em simplicial forest} if every nonempty subcollection
has a leaf, and a {\em simplicial tree} if it is also
connected~\cite{Faridi}.
Every subcollection of a simplicial forest is a simplicial forest.
A vertex is {\em free} if it belongs to exactly one facet.
If $G$ is a joint of a leaf $F$, then every vertex of $F\setminus G$
is free.

Let $\Delta$ be pure of dimension $d\geq1$.
Two facets are {\em adjacent in codimension one}, or
{\em codimension one neighbours}, if they share $d$ vertices.
A {\em proper chain} from $F$ to $G$ is a sequence of facets
$F=F_0,F_1,\ldots,F_s=G$ in which consecutive facets are adjacent in
codimension one.
The complex is {\em connected in codimension one} if any two facets are
joined by a proper chain. Its {\em components in codimension one} are
the equivalence classes of facets under this relation, regarded as
subcollections.
We say that $F$ and $G$ are {\em connected in codimension one inside
$\Delta_{F\cup G}$} if there is a proper chain from $F$ to $G$ all of
whose facets are contained in $F\cup G$.

A {\em $k$-matching} of $\Delta$ is a set of $k$ pairwise disjoint facets.
The {\em matching number} $\nu(\Delta)$ is the largest cardinality of
a matching. We regard the empty set as the unique $0$-matching.
For a matching $\mathcal M$, set
$V_{\mathcal M}=\bigcup_{F\in\mathcal M}F$ and
$x_{\mathcal M}=x_{V_{\mathcal M}}$.
Two disjoint facets $F,G$ form a {\em gap} if
$\Delta_{F\cup G}=\langle F,G\rangle$.
A matching is {\em restricted} if one of its facets forms a gap with
every other facet in the matching. The {\em restricted matching number}
$\nu_0(\Delta)$ is the largest cardinality of a restricted matching.
The {\em $k$-th squarefree power} of $I(\Delta)$ is
\[
I(\Delta)^{[k]}=(x_{\mathcal M}:\mathcal M\text{ is a }k\text{-matching}).
\]
We consider only nonzero squarefree powers of positive order, so
$1\leq p\leq\nu(\Delta)$, and use the convention $I(\Delta)^{[0]}=S$
in inductive arguments.


\begin{Theorem}\label{thm:KNQ}
Let $\Delta$ be a simplicial forest.
\begin{enumerate}
\item[\textnormal{(a)}] If $\mathcal M,\mathcal N$ are $k$-matchings with $V_{\mathcal M}=V_{\mathcal N}$, then $\mathcal M=\mathcal N$. Hence the monomials $x_{\mathcal M}$ are the minimal generators of $I(\Delta)^{[k]}$ \textnormal{\cite[Proposition~1.2, Corollary~1.3]{KNQ24squarefree}}.
\item[\textnormal{(b)}] $\nu(\Delta)-\nu_0(\Delta)\leq 1$ \textnormal{\cite[Proposition~1.4]{KNQ24squarefree}}.
\item[\textnormal{(c)}] If $\Delta$ is pure of positive dimension, then $I(\Delta)^{[k]}$ is not linearly related for $1\leq k<\nu_0(\Delta)$ \textnormal{\cite[Theorem~3.3]{KNQ24squarefree}}.
\end{enumerate}
\end{Theorem}


Let $I$ be a monomial ideal generated in degree $D$, and let $G(I)$
denote its minimal monomial generating set.
We write $\beta_{i,j}(I)$ for its {\em graded Betti numbers}.
The ideal has a {\em linear resolution} if $\beta_{i,j}(I)=0$ whenever
$j\ne D+i$, and it is {\em linearly related} if $\beta_{1,j}(I)=0$
whenever $j\ne D+1$.
Let $G_I$ be the graph on $G(I)$ in which distinct vertices $u,v$ are
adjacent if $\deg\lcm(u,v)=D+1$, and let $G_I^{(u,v)}$ be its induced
subgraph on the generators dividing $\lcm(u,v)$.


\begin{Theorem}[{\cite[Corollary~2.2]{BHZ18powersedgeideals}}]\label{thm:linearlyrelated}
The ideal $I$ is linearly related if and only if, for all $u,v\in G(I)$, there is a path from $u$ to $v$ in $G_I^{(u,v)}$.
\end{Theorem}

Let $\Delta$ be a pure simplicial forest of dimension $d\geq 1$ on the vertex set $V$, and $p\geq 1$. Two $p$-matchings $\mathcal M,\mathcal N$ are \emph{adjacent} if $|V_{\mathcal M}\setminus V_{\mathcal N}|=|V_{\mathcal N}\setminus V_{\mathcal M}|=1$, and for $W\subseteq V$ we denote by $G_p(\Delta;W)$ the graph on the $p$-matchings $\mathcal N$ with $V_{\mathcal N}\subseteq W$, with this adjacency. By Theorem~\ref{thm:KNQ}(a), $x_{\mathcal M}$ and $x_{\mathcal N}$ are adjacent in the graph of Theorem~\ref{thm:linearlyrelated} exactly when $\mathcal M$ and $\mathcal N$ are adjacent, and $\lcm(x_{\mathcal M},x_{\mathcal N})=x_{V_{\mathcal M}\cup V_{\mathcal N}}$. Therefore $I(\Delta)^{[p]}$ is linearly related if and only if
\begin{equation}\label{eq:linrel}
\text{any two $p$-matchings $\mathcal M,\mathcal N$ are connected in }G_p(\Delta;V_{\mathcal M}\cup V_{\mathcal N}).
\end{equation}


For matchings $\mathcal M,\mathcal N$ of a simplicial complex, their
{\em bipartite intersection graph} $B(\mathcal M,\mathcal N)$ has vertex
set $(\mathcal M\setminus\mathcal N)\sqcup(\mathcal N\setminus\mathcal M)$,
with an edge between two facets if they intersect.
The next lemma extends~\cite[Lemma~3.3]{NQ}.


\begin{Lemma}\label{lem:bipartite-intersection-forest}
If $\Delta$ is a simplicial forest and $\mathcal M,\mathcal N$ are matchings of $\Delta$, then $B(\mathcal M,\mathcal N)$ is a forest.
\end{Lemma}

\begin{proof}
Let $P_1,Q_1,\ldots,P_r,Q_r$ be a cycle, with $r\geq 2$, $P_i\in\mathcal M$ and $Q_i\in\mathcal N$, and let $\Gamma=\langle P_1,\ldots,P_r,Q_1,\ldots,Q_r\rangle$. A joint of $P_i$ in $\Gamma$ would contain the nonempty sets $P_i\cap Q_{i-1}$ and $P_i\cap Q_i$ (indices modulo $r$). This is impossible, since the other facets of $\mathcal M$ are disjoint from $P_i$ and every facet of $\mathcal N$ misses $Q_{i-1}$ or $Q_i$. By symmetry no $Q_i$ is a leaf of $\Gamma$, a contradiction.
\end{proof}


The next three lemmas will be used to prove the necessity of the combinatorial conditions in Theorem~\ref{thm:lin-rel-char}.
\begin{Lemma}\label{lem:adjacent}
Let $\mathcal M,\mathcal N$ be adjacent $p$-matchings of $\Delta$, with $V_{\mathcal M}\setminus V_{\mathcal N}=\{a\}$ and $V_{\mathcal N}\setminus V_{\mathcal M}=\{b\}$. Then $B(\mathcal M,\mathcal N)$ is a path $P_1,Q_1,P_2,\ldots,P_r,Q_r$ with $a\in P_1$, $b\in Q_r$, $|P_i\cap Q_i|=d$ for all $i$ and $|Q_i\cap P_{i+1}|=1$ for $i<r$. In particular, if $\Gamma$ is an induced subcollection containing $\mathcal M$ and $\mathcal N$, then $\mathcal M$ and $\mathcal N$ have the same number of facets in each component of $\Gamma$ in codimension one.
\end{Lemma}

\begin{proof}
A facet $P\in\mathcal M\setminus\mathcal N$ not containing $a$ is contained in $V_{\mathcal N}$ and is disjoint from the facets of $\mathcal M\cap\mathcal N$, so it is covered by $\mathcal N\setminus\mathcal M$; as $P\notin\mathcal N$ and $\Delta$ is pure, it meets at least two facets of $\mathcal N\setminus\mathcal M$. The same holds for the facets of $\mathcal N\setminus\mathcal M$ not containing $b$, while the facets containing $a$ and $b$ have degree at least one because $d\geq 1$. Hence the forest $B(\mathcal M,\mathcal N)$ has no isolated vertices and at most two vertices of degree one, so it is a path whose ends are the facets $P_1\ni a$ and $Q_r\ni b$. Now $P_1\setminus\{a\}\subseteq Q_1$, $Q_i\subseteq P_i\cup P_{i+1}$ for $i<r$ and $P_i\subseteq Q_{i-1}\cup Q_i$ for $i>1$, which gives the stated cardinalities. The last assertion follows, since $P_i$ and $Q_i$ lie in the same component of $\Gamma$ in codimension one.
\end{proof}


\begin{Lemma}\label{lem:unique-facet-overlap}
Let $\mathcal A$ be a $(p-1)$-matching and let $F,G$ be facets disjoint from $V_{\mathcal A}$ with $F\cap G\neq\emptyset$. Then every $p$-matching $\mathcal M$ with $V_{\mathcal M}\subseteq V_{\mathcal A}\cup F\cup G$ contains exactly one facet disjoint from $V_{\mathcal A}$, and this facet lies in $F\cup G$.
\end{Lemma}

\begin{proof}
Two disjoint facets cannot fit in $F\cup G$, since $|F\cup G|<2(d+1)$; this gives uniqueness. Suppose that every facet of $\mathcal M$ meets $V_{\mathcal A}$. Then $\mathcal M\setminus\mathcal A$ has one facet more than $\mathcal A\setminus\mathcal M$, and each of its facets meets a facet of $\mathcal A\setminus\mathcal M$. Thus some component of the forest $B(\mathcal M,\mathcal A)$ has more vertices in $\mathcal M$ than in $\mathcal A$, and a degree count shows that it has two leaves $P,Q\in\mathcal M$. The facet $P$ meets only one facet of $\mathcal A$, so by purity it has a vertex outside $V_{\mathcal A}$, hence in $F\cup G$; the same holds for $Q$. If both $P$ and $Q$ met $F$, then $F$, which is not in $\mathcal M$, and the path from $P$ to $Q$ would form a cycle in $B(\mathcal M,\mathcal A\cup\{F\})$, against Lemma~\ref{lem:bipartite-intersection-forest}; the same holds for $G$. So we may assume that $P$ meets $F$ but not $G$, and $Q$ meets $G$ but not $F$. Let $P=P_0,E_1,P_1,\ldots,E_t,P_t=Q$ be the path in $B(\mathcal M,\mathcal A)$ and $\Gamma=\langle P_0,\ldots,P_t,E_1,\ldots,E_t,F,G\rangle$. Each facet of $\Gamma$ meets two facets of $\Gamma$ (namely $E_i$ meets $P_{i-1},P_i$; $P_i$ meets $E_i,E_{i+1}$ for $0<i<t$; $P$ meets $E_1,F$; $Q$ meets $E_t,G$; $F$ meets $P,G$; $G$ meets $Q,F$), and since the facets of $\mathcal M$ are pairwise disjoint, those of $\mathcal A$ as well, $F,G$ avoid $V_{\mathcal A}$, $P\cap G=\emptyset$ and $Q\cap F=\emptyset$, no third facet of $\Gamma$ contains both intersections. So $\Gamma$ has no leaf, a contradiction.
\end{proof}


\begin{Lemma}\label{lem:path-gives-codim-one}
Let $\mathcal A,F,G$ be as in Lemma~\ref{lem:unique-facet-overlap} and set $L=V_{\mathcal A}\cup F\cup G$. If $\mathcal A\cup\{F\}$ and $\mathcal A\cup\{G\}$ are connected in $G_p(\Delta;L)$, then $F$ and $G$ are connected in codimension one inside $\Delta_{F\cup G}$.
\end{Lemma}

\begin{proof}
For a vertex $\mathcal M$ of $G_p(\Delta;L)$ let $H_{\mathcal M}\subseteq F\cup G$ be the facet of $\mathcal M$ disjoint from $V_{\mathcal A}$ given by Lemma~\ref{lem:unique-facet-overlap}. It suffices to show that, if $\mathcal M,\mathcal N$ are adjacent and $H_{\mathcal M}\neq H_{\mathcal N}$, then $|H_{\mathcal M}\cap H_{\mathcal N}|=d$; following a path from $\mathcal A\cup\{F\}$ to $\mathcal A\cup\{G\}$ then gives a proper chain from $F$ to $G$ inside $F\cup G$. By uniqueness, neither $H_{\mathcal M}$ nor $H_{\mathcal N}$ lies in $\mathcal M\cap\mathcal N$, and they intersect, because they are contained in $F\cup G$. So they are adjacent in the path of Lemma~\ref{lem:adjacent}. If they are $P_i,Q_i$ we are done. Otherwise $H_{\mathcal N}=Q_i$ and $H_{\mathcal M}=P_{i+1}$ for some $i<r$, and $(\mathcal M\cap\mathcal N)\cup\{P_1,\ldots,P_i,Q_{i+1},\ldots,Q_r\}$ is a $p$-matching inside $L$ all of whose facets meet $V_{\mathcal A}$, against Lemma~\ref{lem:unique-facet-overlap}.
\end{proof}


\section{Linearly related squarefree powers}\label{sec:linear-related}

Throughout this section, let $\Delta$ be a pure simplicial forest of
dimension $d\geq1$. We characterize the linearly related squarefree
powers of $I(\Delta)$ by two combinatorial conditions.

\begin{Definition}\label{def:conditions}
We say that $\Delta$ satisfies
\begin{enumerate}
\item[\textnormal{(C1)$_p$}] if $\nu(\Delta_{V\setminus(F\cup G)})<p-1$ for all facets $F,G$ with $0<|F\cap G|<d$ that are not connected in codimension one inside $\Delta_{F\cup G}$;
\item[\textnormal{(C2)$_p$}] if, for every $(p+1)$-matching $\mathcal L$ and all $E,E'\in\mathcal L$, the matchings $\mathcal L\setminus\{E\}$ and $\mathcal L\setminus\{E'\}$ are connected in $G_p(\Delta;V_{\mathcal L})$.
\end{enumerate}
\end{Definition}

Condition $(\mathrm C2)_p$ is the instance of~\eqref{eq:linrel}
corresponding to the pairs
$\mathcal L\setminus\{E\},\mathcal L\setminus\{E'\}$, and it is
vacuous when $p=\nu(\Delta)$.
Both conditions are inherited by every induced subcollection
$\Delta_W$, because $(\Delta_W)_{F\cup G}=\Delta_{F\cup G}$ for facets
$F,G$ of $\Delta_W$, and
$G_p(\Delta_W;U)=G_p(\Delta;U)$ for $U\subseteq W$.

The following theorem establishes the equivalence between these two
conditions and the connectivity criterion in~\eqref{eq:linrel}.


\begin{Theorem}\label{thm:lin-rel-char}
Let $1\leq p\leq\nu(\Delta)$. Then $I(\Delta)^{[p]}$ is linearly related if and only if $\Delta$ satisfies \textnormal{(C1)$_p$} and \textnormal{(C2)$_p$}.
\end{Theorem}
\begin{proof}
If $I(\Delta)^{[p]}$ is linearly related, (C2)$_p$ is a special case of \eqref{eq:linrel}, and (C1)$_p$ follows from \eqref{eq:linrel} and Lemma~\ref{lem:path-gives-codim-one}, applied to a $(p-1)$-matching $\mathcal A$ of $\Delta_{V\setminus(F\cup G)}$.

Conversely, let $\mathcal A,\mathcal B$ be $p$-matchings; we show by induction on $h=|\mathcal A\setminus\mathcal B|$ that they are connected in $G_p(\Delta;V_{\mathcal A}\cup V_{\mathcal B})$. Let $h>0$. The forest $B(\mathcal A,\mathcal B)$ has $2h$ vertices, hence at most $2h-1$ edges, so some $E\in\mathcal A\setminus\mathcal B$ meets at most one facet of $\mathcal B$. If $E$ meets exactly one facet $H$ of $\mathcal B$, then the other facets of $\mathcal B$ avoid $E\cup H$, and either $|E\cap H|=d$ or (C1)$_p$ gives a proper chain $H=H_0,\ldots,H_s=E$ inside $\Delta_{E\cup H}$; replacing $H$ by $H_1,\ldots,H_s$ in turn gives a path from $\mathcal B$ to $\mathcal B'=(\mathcal B\setminus\{H\})\cup\{E\}$. If $E$ meets no facet of $\mathcal B$, then $\mathcal B\cup\{E\}$ is a $(p+1)$-matching, and for any $H\in\mathcal B\setminus\mathcal A$ condition (C2)$_p$ joins $\mathcal B$ to $\mathcal B'=(\mathcal B\setminus\{H\})\cup\{E\}$. In both cases the path stays inside $V_{\mathcal A}\cup V_{\mathcal B}$ and $|\mathcal A\setminus\mathcal B'|=h-1$, so we conclude by induction and \eqref{eq:linrel}.
\end{proof}


\begin{Remark}\label{rem:consequences}
If $\Delta$ satisfies (C1)$_p$ and (C2)$_p$, then $\nu(\Delta)-1\leq\nu_0(\Delta)\leq p\leq\nu(\Delta)$ by Theorem~\ref{thm:lin-rel-char} and Theorem~\ref{thm:KNQ}(b),(c). So the conditions can only hold for $p\in\{\nu(\Delta)-1,\nu(\Delta)\}$.
\end{Remark}


\begin{Proposition}\label{prop:C2}
If $\Delta$ satisfies \textnormal{(C2)$_p$}, then $\Delta_{V_{\mathcal L}}$ is connected in codimension one for every $(p+1)$-matching $\mathcal L$. In particular, this holds when $I(\Delta)^{[p]}$ is linearly related.
\end{Proposition}

\begin{proof}
Let $\Gamma=\Delta_{V_{\mathcal L}}$ and $E,E'\in\mathcal L$. All the matchings along a path from $\mathcal L\setminus\{E\}$ to $\mathcal L\setminus\{E'\}$ in $G_p(\Delta;V_{\mathcal L})$ are matchings of $\Gamma$, so by Lemma~\ref{lem:adjacent} they have the same number of facets in each component of $\Gamma$ in codimension one. Hence $E$ and $E'$ lie in the same component $\mathcal C$, and $\mathcal C$ covers $V_{\mathcal L}$. Let $n=|V_{\mathcal L}|$. Ordering the facets of $\mathcal C$ so that each of them is adjacent in codimension one to an earlier one, each facet adds at most one new vertex, so $\mathcal C$ has at least $n-d$ facets. On the other hand, a pure simplicial forest with $m\geq 1$ facets of size $d+1$ on $n$ vertices satisfies $m\leq n-d$, by induction on $m$: removing a leaf removes at least one free vertex. Hence $\Gamma=\mathcal C$.
\end{proof}


\begin{Remark}\label{rem:NQT}
Under \textnormal{(C1)$_p$} and \textnormal{(C2)$_p$},
Proposition~\ref{prop:C2} is nonvacuous only when
$p=\nu(\Delta)-1=\nu_0(\Delta)$, by Remark~\ref{rem:consequences}.
When $d=1$, the complex $\Delta$ is a graph-theoretic forest, and
the conclusion asserts that every $(p+1)$-matching covers a connected
induced subgraph. The same connectivity condition characterizes
linearly related squarefree powers of edge ideals of arbitrary graphs
$G$, for $2\leq p\leq\nu(G)$
\cite[Theorem~3.1]{NQT26}.
\end{Remark}


\section{From linear relations to linear quotients}\label{sec:LQ}

Throughout this section, let $\Delta$ be a pure simplicial forest of dimension $d\geq1$. 

A monomial ideal $I$ has {\em linear quotients} if its minimal monomial
generators admit an order $u_1,\ldots,u_m$ such that
$(u_1,\ldots,u_{j-1}):u_j$ is generated by variables for every $j\geq2$.
Such an order is called a {\em linear quotient order}.
Equivalently, for every $i<j$ there are $k<j$ and a variable $x$ such that $u_k/\gcd(u_k,u_j)=x$ and $x$ divides $u_i/\gcd(u_i,u_j)$; see~\cite[Lemma~8.2.3]{HHmonomialideals}.
For an ideal generated in one degree, linear quotients imply a linear
resolution, and a linear resolution implies that the ideal is linearly
related. 
In view of Theorem~\ref{thm:lin-rel-char}, it therefore remains
to prove that $(\mathrm C1)_p$ and $(\mathrm C2)_p$ imply the existence
of a linear quotient order for $I(\Delta)^{[p]}$.


\subsection{Good leaves and end leaves}\label{subsec:good-leaves}

A facet $F$ of $\Delta$ is a {\em good leaf} if the sets
$F\cap G$, $G\in\mathcal F(\Delta)$, are totally ordered by inclusion.
A good leaf remains a good leaf in every subcollection containing it,
and every nonempty simplicial forest has a good
leaf~\cite[Corollary~3.4]{HHTZ}. A proper chain from $E$ to $H$ is {\em irredundant} if no proper
subsequence is a proper chain from $E$ to $H$.


\begin{Lemma}\label{lem:chain}
Let $E=H_0,H_1,\ldots,H_s=H$ be an irredundant proper chain in $\Delta_{E\cup H}$, and write $H_{i+1}=(H_i\setminus\{a_i\})\cup\{b_i\}$. Then $a_i\in E\setminus H$ and $b_i\in H\setminus E$ for every $i$, and $E\cap H\subseteq H_i$ for every $i$.
\end{Lemma}
\begin{proof}
The tree $\langle H_i,\ldots,H_s\rangle$ is connected in codimension one, so it has two leaves \cite[Proposition~1.16]{Zheng2004}, while $H_{i+1},\ldots,H_{s-1}$ are not leaves since the chain is irredundant \cite[Lemma~1.10]{Zheng2004}. So $H_i$ is a leaf, and a joint $J$ satisfies $H_i\cap J=H_i\cap H_{i+1}$; hence $a_i$ lies in none of $H_{i+1},\ldots,H_s$, and $a_i\in E\setminus H$ because $H_i\subseteq E\cup H$. Symmetrically $b_i\in H\setminus E$. As $a_i\notin H$, induction gives $E\cap H\subseteq H_i$.
\end{proof}


A {\em special cycle} of length $s\geq3$ is a sequence
$F_1,v_1,\ldots,F_s,v_s$ of distinct facets and distinct vertices with $v_i\in F_i\cap F_{i+1}$, where $F_{s+1}=F_1$, such that none of $F_1,\ldots,F_s$ contains more than two of $v_1,\ldots,v_s$.
A simplicial complex is a simplicial forest if and only if it has no special cycle of length at least three~\cite[Theorem~3.2]{HHTZ}.

\begin{Lemma}\label{lem:two-good-leaves}
Every simplicial forest $\Delta$ with at least two facets has two good leaves.
\end{Lemma}

\begin{proof}
Let $G$ be a good leaf of $\Delta$. If $G$ is disjoint from the other facets, a good leaf of $\langle\mathcal F(\Delta)\setminus\{G\}\rangle$ is a good leaf of $\Delta$. Otherwise let $A$ be the largest of the sets $G\cap H$, $H\neq G$, attained at $H=K$, and choose $g\in G\setminus A$; then $g$ is free and $|G|\geq 2$. Let $\Delta_1,\Delta_2$ be two copies of $\Delta$ sharing only the vertex $g$, with $H_1,H_2$ the copies of a facet $H$, and let $\Delta'$ be the complex with facets $\mathcal F(\Delta_1)\cup\mathcal F(\Delta_2)$. Only $G_1$ and $G_2$ meet across the two copies, in $\{g\}$, so a cycle of $\Delta'$ using both copies would pass from one to the other at least twice, each time between $G_1$ and $G_2$, which is impossible for a cycle of length at least three. A special cycle inside one copy is excluded since $\Delta$ is a forest, so $\Delta'$ is a forest by \cite[Theorem~3.2]{HHTZ}. Since $G_i\cap G_{3-i}=\{g\}$ and $G_i\cap K_i=A_i$ are incomparable, the good leaf of $\Delta'$ given by \cite[Corollary~3.4]{HHTZ} is $L_i$ for some facet $L\neq G$. As $g\notin L$, the facet $L_i$ misses the other copy, so $L$ is a good leaf of $\Delta$.
\end{proof}


Let $F$ be a good leaf of $\Delta$. Two facets meeting $F$ meet each other, so a matching contains at most one facet meeting $F$, and $V_{\mathcal M}\cap F$ is either empty or of the form $G\cap F$. A \emph{ridge} of $\Delta$ is a set of $d$ vertices contained in at least two facets. If $F$ has a codimension one neighbour, the sets $F\cap G$ with $|F\cap G|=d$ are comparable of the same size, hence equal to one ridge $C$; we write $F=C\cup\{f\}$, where $f$ is free. The \emph{clones} of $F$ are the facets $C\cup\{x\}$ with $x$ free; if $F$ has no codimension one neighbour, $F$ is its only clone. A clone $F'$ satisfies $F'\cap H=F\cap H$ for every facet $H$ which is not a clone, and in particular $\Delta_{V\setminus F'}=\Delta_{V\setminus F}$.

\begin{Definition}\label{def:end-leaf}
A good leaf $F$ is an \emph{end leaf} if at most one of its codimension one neighbours is not a clone of $F$.
\end{Definition}



To locate an end leaf, we use the incidence graph of facets and ridges. Let $\mathcal R(\Delta)$ be the bipartite graph on the facets and the ridges of $\Delta$, in which a facet is adjacent to the ridges it contains. Two facets are adjacent in codimension one if and only if they have a common neighbour in $\mathcal R(\Delta)$.

\begin{Lemma}\label{lem:ridge-graph}
\begin{enumerate}
\item $\mathcal R(\Delta)$ is a forest.
\item If $E_0,C_1,E_1,\ldots,C_k,E_k$ is a path in $\mathcal R(\Delta)$, with $E_i$ facets, then $E_0\cap E_k\subseteq E_i$ for every $i$.
\item If $\Delta$ is connected in codimension one and $H$ is a facet containing only one ridge $C$, then the vertex of $H\setminus C$ is free.
\end{enumerate}
\end{Lemma}

\begin{proof}
(1) In a cycle $F_1,C_1,\ldots,F_k,C_k$ each facet is $F_i=C_{i-1}\cup C_i$; for $k=2$ this gives $F_1=F_2$, and for $k\geq 3$ no $F_i$ is a leaf of $\langle F_1,\ldots,F_k\rangle$. (2) By induction on $k$. For $0<i<k$ the facet $E_i$ is not a leaf of $\langle E_0,\ldots,E_k\rangle$, as in (1), so we may assume that $E_k$ is a leaf; its joint contains $C_k$, so by (1) it is $E_{k-1}$. Thus $E_0\cap E_k\subseteq E_0\cap E_{k-1}$, and we apply induction to $E_0,\ldots,E_{k-1}$. (3) If $h\in H\setminus C$ lies in a facet $K\neq H$, the path from $H$ to $K$ in $\mathcal R(\Delta)$ starts with $H,C,X$ for a facet $X\supseteq C$, $X\neq H$, and $X\neq K$ because $C\not\subseteq K$. By (2), $h\in X$, so $H\subseteq X$, a contradiction.
\end{proof}


The preceding properties of the ridge graph allow us to prove the existence of end leaves.
\begin{Lemma}\label{lem:end-leaf}
If $\Delta$ has at least two facets and is connected in codimension one, then it has an end leaf.
\end{Lemma}
\begin{proof}
Every good leaf $F$ contains exactly one ridge $C_F$, and by Lemma~\ref{lem:ridge-graph}(3) the facets whose only ridge is $C_F$ are clones of $F$. So we need a good leaf $F$ such that at most one facet containing $C_F$ contains a second ridge. Fix a vertex $\rho$ of $\mathcal R(\Delta)$ and choose a good leaf $F$ with $C_F$ as far as possible from $\rho$. Suppose that two facets containing $C_F$ contain a second ridge, and let $G_1$ be one of them not lying on the path from $C_F$ to $\rho$, with a second ridge $C'$. Let $\mathcal S$ be the set of facets in the component of $\mathcal R(\Delta)\setminus\{C_F\}$ containing $G_1$, and $\Sigma=\langle\mathcal S\cup\{F\}\rangle$. All the vertices of this component are farther from $\rho$ than $C_F$, and $G_1$ is the only facet of $\mathcal S$ containing $C_F$. Since $G_1$ is not a leaf of $\Sigma$ (a joint would contain $C_F\cup C'=G_1$), Lemma~\ref{lem:two-good-leaves} gives a good leaf $H\notin\{F,G_1\}$ of $\Sigma$.

We claim that $H$ is a good leaf of $\Delta$. For a facet $K\notin\Sigma$, the path from $H$ to $K$ in $\mathcal R(\Delta)$ passes through $G_1$, $C_F$ and a facet $X\supseteq C_F$ different from $G_1$, so $H\cap K\subseteq G_1\cap X=C_F\subseteq F$ by Lemma~\ref{lem:ridge-graph}(2), that is, $H\cap K=D\cap(F\cap K)$ with $D=H\cap F$. These sets form a chain because $F$ is a good leaf. If $P\in\Sigma$, then $H\cap P$ either contains $D$, and hence all these sets, or is contained in $D$, and then $H\cap P=D\cap(F\cap P)$ is comparable with them. This proves the claim. But $C_H\neq C_F$ lies in the component of $G_1$, so it is farther from $\rho$ than $C_F$, against the choice of $F$.
\end{proof}


\subsection{Linear quotient and the main result}\label{subsec:lq-construction}

We first establish a reduction that allows the inductive argument to
be restricted to facets occurring in $p$-matchings.

A facet is {\em $p$-active} if it belongs to a $p$-matching.
Deleting all other facets preserves $I(\Delta)^{[p]}$ and
$\nu(\Delta)$, since every facet of a maximum matching is $p$-active.
It also preserves \textnormal{(C1)$_p$} and \textnormal{(C2)$_p$}:
the $p$-matching graphs are unchanged, and every facet contained in
$F\cup G$ is $p$-active whenever a $(p-1)$-matching avoids $F\cup G$.
Hence the proper chains required by \textnormal{(C1)$_p$} are retained.


\begin{Lemma}\label{lem:two-facets}
Let $\mathcal L$ be a $(p+1)$-matching and $\mathcal N$ a $p$-matching with $V_{\mathcal N}\subseteq V_{\mathcal L}$. Then every facet of $\mathcal N\setminus\mathcal L$ meets exactly two facets of $\mathcal L$ and is contained in their union.
\end{Lemma}
\begin{proof}
Let $s=|\mathcal L\setminus\mathcal N|$. A facet of $\mathcal N\setminus\mathcal L$ is covered by $\mathcal L\setminus\mathcal N$ and, by purity, meets at least two of its facets. The forest $B(\mathcal L,\mathcal N)$ has $2s-1$ vertices, hence at most $2s-2$ edges, and each edge has one end in $\mathcal N\setminus\mathcal L$, which has $s-1$ elements. So all these degrees are exactly two.
\end{proof}


\begin{Lemma}\label{lem:bridge}
Let $F,P$ be disjoint facets and let $\mathcal B$ be the set of facets contained in $F\cup P$. Then the sets $K\cap F$, $K\in\mathcal B$, are totally ordered, and $K$ is determined by $K\cap F$. In particular, if $K'\cap F=(K\cap F)\cup\{b\}$ with $b\notin K$, then $K'=(K\setminus\{a\})\cup\{b\}$ for some $a\in K\cap P$.
\end{Lemma}

\begin{proof}
Write $K=A_K\sqcup B_K$ with $A_K=K\cap F$ and $B_K=K\cap P$. Suppose that $A_K$ and $A_{K'}$ are incomparable; then all four sets $A_K,B_K,A_{K'},B_{K'}$ are nonempty. In $\langle F,P,K,K'\rangle$ neither $F$ nor $K$ nor $K'$ is a leaf: a joint of $F$ would contain $A_K\cup A_{K'}$, and a joint of $K$ would give $B_K\subseteq A_K$, $A_K\subseteq B_K$ or $K\subseteq K'$, so $P$ is a leaf, and we may assume $B_{K'}\subseteq B_K$. Then $K\cap K'\supseteq B_{K'}\neq\emptyset$, and one checks in the same way that $\langle F,K,K'\rangle$ has no leaf, a contradiction. So the $A_K$ are totally ordered, and by symmetry so are the $B_K$, in the opposite order since $|A_K|+|B_K|=d+1$. This gives the remaining assertions.
\end{proof}


\begin{Lemma}\label{lem:choice}
Assume that $\Delta$ has at least two facets, all of them $p$-active, that $\nu(\Delta)\in\{p,p+1\}$, and that $\Delta$ satisfies \textnormal{(C1)$_p$} and \textnormal{(C2)$_p$}. Then there are a good leaf $F$ of $\Delta$ and a set $W\subseteq V\setminus F$ with $|W|\leq 1$ such that, for $\Gamma=\Delta_{V\setminus F}$ and $\Gamma'=\Delta_{V\setminus(F\cup W)}$:
\begin{enumerate}
\item[\textnormal{(i)}] every facet which meets $F$ and is not a clone of $F$ contains $W$;
\item[\textnormal{(ii)}] $\nu(\Gamma')\leq p-1$;
\item[\textnormal{(iii)}] if $p\geq 2$, then $\Gamma$ satisfies \textnormal{(C1)$_{p-1}$} and \textnormal{(C2)$_{p-1}$}.
\end{enumerate}
If $\nu(\Delta)=p$ one can take any good leaf $F$ and $W=\emptyset$; if $\nu(\Delta)=p+1$ one can take any end leaf $F$ and $W=\{g\}$, where $C\cup\{g\}$ is the unique codimension one neighbour of $F$ which is not a clone.
\end{Lemma}

\begin{proof}
Suppose that $p\geq2$. For any facet $F$, the induced subcollection
$\Gamma$ satisfies $(\mathrm C1)_{p-1}$.
Indeed, if $\mathcal S$ is a $(p-2)$-matching of $\Gamma$ disjoint from
$E\cup H$, then $\mathcal S\cup\{F\}$ is a $(p-1)$-matching of
$\Delta$ disjoint from $E\cup H$, and $(\mathrm C1)_p$ applies. If $F$ is a good leaf, then
\begin{equation}\label{eq:good-leaf-matching-drop}
\nu(\Gamma)=\nu(\Delta)-1,
\end{equation}
since $F$ can be added to every matching of $\Gamma$, while a maximum matching of $\Delta$ contains exactly one facet meeting $F$, which can be replaced by $F$.

If $\nu(\Delta)=p$, choose any good leaf $F$ and set $W=\varnothing$.
Then (i) holds trivially, and (ii) follows
from~\eqref{eq:good-leaf-matching-drop}.
When $p\geq2$, condition $(\mathrm C2)_{p-1}$ for $\Gamma$ is vacuous
by~\eqref{eq:good-leaf-matching-drop}. Together with the preceding
verification of $(\mathrm C1)_{p-1}$, this proves (iii).

Let $\nu(\Delta)=p+1$. First, $\Delta$ is connected in codimension one. Indeed, by Theorem~\ref{thm:lin-rel-char} and \eqref{eq:linrel}, any two $p$-matchings are connected in $G_p(\Delta;V)$, so by Lemma~\ref{lem:adjacent} all $p$-matchings have the same number of facets in each component of $\Delta$ in codimension one. Applying this to $\mathcal L\setminus\{E\}$, $E\in\mathcal L$, for a $(p+1)$-matching $\mathcal L$, all $p$-matchings lie in one component, which then contains every facet since all facets are $p$-active. By Lemma~\ref{lem:end-leaf} there is an end leaf $F=C\cup\{f\}$. Not all its codimension one neighbours are clones, since otherwise every facet would contain $C$ and $\nu(\Delta)=1$; let $G=C\cup\{g\}$ be the non-clone one and $W=\{g\}$.

(i) Let $H$ meet $F$, not a clone. If $|H\cap F|=d$, then $H=G$. Otherwise, the other facets of a $p$-matching containing $H$ avoid $F\cup H$, so (C1)$_p$ and Lemma~\ref{lem:chain} give an irredundant proper chain $F,H_1,\ldots,H$ inside $F\cup H$ with $H_1=C\cup\{y\}$ and $y\in H$. As $y$ is not free, $H_1=G$ and $g\in H$.

(ii) Let $\mathcal R$ be a $p$-matching of $\Gamma$ with $g\notin V_{\mathcal R}$, and $\mathcal L=\mathcal R\cup\{F\}$. By (i), no facet inside $V_{\mathcal L}$ other than $F$ meets $F$ (clones $F'\neq F$ have a free vertex outside $V_{\mathcal L}$). So every vertex of $G_p(\Delta;V_{\mathcal L})$ other than $\mathcal R$ contains $F$ (by Theorem~\ref{thm:KNQ}(a)), and $\mathcal R$ is isolated since $|F|\geq 2$, against (C2)$_p$. Hence $g$ is covered by every $p$-matching of $\Gamma$, and $\nu(\Gamma')\leq p-1$.

(iii) We prove (C2)$_{p-1}$ for $\Gamma$. Let $\mathcal R$ be a $p$-matching of $\Gamma$, $E,E'\in\mathcal R$, let $P\in\mathcal R$ be the facet containing $g$ and $\mathcal L=\mathcal R\cup\{F\}$. For a vertex $\mathcal N$ of $G_p(\Delta;V_{\mathcal L})$ let $H_{\mathcal N}$ be its facet meeting $F$, and $H_{\mathcal R}=P$. We have $H_{\mathcal N}\subseteq F\cup P$: otherwise $H_{\mathcal N}\notin\mathcal L$ is not a clone, so it contains $g$ by (i), and by Lemma~\ref{lem:two-facets} it lies in the union of $F$ and $P$. By Lemma~\ref{lem:bridge}, $H_{\mathcal N}$ is determined by $V_{\mathcal N}\cap F=H_{\mathcal N}\cap F$, and these sets are totally ordered. Put $\overline{\mathcal N}=\mathcal N\setminus\{H_{\mathcal N}\}$, a $(p-1)$-matching of $\Gamma$ inside $V_{\mathcal R}$; note that $\overline{\mathcal L\setminus\{E\}}=\mathcal R\setminus\{E\}$. We show that if $\mathcal N,\mathcal N'$ are adjacent, with $V_{\mathcal N'}=(V_{\mathcal N}\setminus\{c\})\cup\{b\}$, then $\overline{\mathcal N}$ and $\overline{\mathcal N'}$ are equal or adjacent; then (C2)$_p$ gives a walk from $\mathcal R\setminus\{E\}$ to $\mathcal R\setminus\{E'\}$ in $G_{p-1}(\Gamma;V_{\mathcal R})$. If $H_{\mathcal N}=H_{\mathcal N'}$, the claim is clear. Otherwise the traces on $F$ differ, so exactly one of $b,c$ lies in $F$ (two comparable traces of the same size coincide), say $b$, and $H_{\mathcal N'}=(H_{\mathcal N}\setminus\{a\})\cup\{b\}$ with $a\in H_{\mathcal N}\cap P$ by Lemma~\ref{lem:bridge}. If $c\in H_{\mathcal N}$, then $c=a$ and $\overline{\mathcal N'}=\overline{\mathcal N}$; if not, $V_{\overline{\mathcal N'}}=(V_{\overline{\mathcal N}}\setminus\{c\})\cup\{a\}$.
\end{proof}


We now use Lemma~\ref{lem:choice} to construct a linear quotient order required to complete the proof of the main theorem.

\begin{Theorem}\label{thm:main}
Let $\Delta$ be a pure simplicial forest of dimension $d\geq 1$ and let $1\leq p\leq\nu(\Delta)$. The following are equivalent:
\begin{enumerate}
\item[\textnormal{(a)}] $I(\Delta)^{[p]}$ has linear quotients;
\item[\textnormal{(b)}] $I(\Delta)^{[p]}$ has a linear resolution;
\item[\textnormal{(c)}] $I(\Delta)^{[p]}$ is linearly related;
\item[\textnormal{(d)}] $\Delta$ satisfies \textnormal{(C1)$_p$} and \textnormal{(C2)$_p$}.
\end{enumerate}
In this case $\nu(\Delta)-1\leq\nu_0(\Delta)\leq p$.
\end{Theorem}

\begin{proof}
Since $I(\Delta)^{[p]}$ is generated in degree $p(d+1)$, linear quotients imply a linear resolution, and a linear resolution implies that the ideal is linearly related. Thus (a)$\Rightarrow$(b)$\Rightarrow$(c). The equivalence (c)$\Leftrightarrow$(d) follows from Theorem~\ref{thm:lin-rel-char}, while the numerical assertion follows from Remark~\ref{rem:consequences}. It therefore remains to prove (d)$\Rightarrow$(a). We use induction on the number of facets of $\Delta$, for all $p$ at once, the case of one facet being trivial. By the reduction at the beginning of
Subsection~\ref{subsec:lq-construction}, we may assume that every facet
is $p$-active. Remark~\ref{rem:consequences} then gives
$\nu(\Delta)\in\{p,p+1\}$. Since the generators are squarefree of the same degree, an order of the $p$-matchings gives linear quotients if and only if for all $\mathcal M$ before $\mathcal N$ there is $\mathcal N'$ before $\mathcal N$ with
\begin{equation}\label{eq:lq}
V_{\mathcal N'}=(V_{\mathcal N}\setminus\{a\})\cup\{b\},
\qquad a\in V_{\mathcal N},\quad
b\in V_{\mathcal M}\setminus V_{\mathcal N}.
\end{equation}
We also use that multiplying an ideal with linear quotients by a monomial preserves linear quotients. Let $F$, $W$, $\Gamma$, $\Gamma'$ be as in Lemma~\ref{lem:choice}.

If $F$ has no codimension one neighbour, then no other facet $H$ meets $F$: the other facets of a $p$-matching containing $H$ avoid $F\cup H$, so (C1)$_p$ would give a proper chain from $F$ to $H$. So $\nu(\Delta)=p$ (for $\nu(\Delta)=p+1$, $F$ is an end leaf of a complex connected in codimension one), every $p$-matching contains $F$ by \eqref{eq:good-leaf-matching-drop}, and $I(\Delta)^{[p]}=x_F\,I(\Gamma)^{[p-1]}$ has linear quotients by Lemma~\ref{lem:choice}(iii) and induction.

Otherwise let $F_i=C\cup\{f_i\}$, $1\leq i\leq t$, be the clones of $F$, with $F_1=F$, and let $A_1\subsetneq\cdots\subsetneq A_\ell$ be the distinct nonempty sets $H\cap F$ with $H$ not a clone; they are contained in $C$. We call a facet $H$ which is not a clone and satisfies $H\cap F=A$ an \emph{$A$-facet}. Every $p$-matching lies in exactly one of the blocks $\mathcal Q_i=\{\mathcal M:F_i\in\mathcal M\}$, $\mathcal Q_{A_j}=\{\mathcal M:\mathcal M\text{ contains an }A_j\text{-facet}\}$ and $\mathcal Q_0=\{\mathcal M:V_{\mathcal M}\cap F=\emptyset\}$.

\emph{Each block has linear quotients.} The ideal of $\mathcal Q_i$ is $x_{F_i}I(\Gamma)^{[p-1]}$, and we use Lemma~\ref{lem:choice}(iii) and induction (it is principal if $p=1$). The block $\mathcal Q_0$ is empty if $\nu(\Delta)=p$, and it consists of the $p$-matchings of $\Gamma$ if $\nu(\Delta)=p+1$, which satisfies (C1)$_p$ and (C2)$_p$; so induction applies. Fix $A=A_j$ and let $\Gamma_A$ be the subcollection whose facets are the $A$-facets and the facets of $\Gamma'$. By Lemma~\ref{lem:choice}(i) the $A$-facets contain $W$, so the other facets of a matching in $\mathcal Q_A$ lie in $\Gamma'$; conversely, a matching of $\Gamma_A$ has at most one $A$-facet and, by Lemma~\ref{lem:choice}(ii), at most $p-1$ facets of $\Gamma'$. Hence $\mathcal Q_A$, which is nonempty because the $A$-facets are $p$-active, is the set of $p$-matchings of $\Gamma_A$, so $\nu(\Gamma_A)=p$ and (C2)$_p$ is empty for $\Gamma_A$. For (C1)$_p$, let $P,Q$ be facets of $\Gamma_A$ with $0<|P\cap Q|<d$ and $\mathcal S$ a $(p-1)$-matching of $\Gamma_A$ avoiding $P\cup Q$. If $P$ is an $A$-facet and $Q\in\Gamma'$, then $\mathcal S\cup\{Q\}$ would be a $p$-matching of $\Gamma'$. Otherwise (C1)$_p$ for $\Delta$ gives an irredundant proper chain from $P$ to $Q$ inside $P\cup Q$. If $P,Q\in\Gamma'$, it lies in $\Gamma'$. If both are $A$-facets, then by Lemma~\ref{lem:chain} each facet of the chain contains $A$, meets $F$ in $A$ and is not a clone, so it is an $A$-facet. Thus induction applies to $\Gamma_A$, which does not contain $F$.

\emph{Gluing the blocks.} We order the matchings block by block, in the order
\[
\mathcal Q_1,\ldots,\mathcal Q_t,\mathcal Q_{A_\ell},\ldots,\mathcal Q_{A_1},\mathcal Q_0,
\]
and check \eqref{eq:lq} for $\mathcal M$ in an earlier block than $\mathcal N$.

If $\mathcal M\in\mathcal Q_i$ and $\mathcal N\in\mathcal Q_j$ with $i<j$, take $\mathcal N'=(\mathcal N\setminus\{F_j\})\cup\{F_i\}$ and $b=f_i$.

Let $\mathcal N\in\mathcal Q_A$, with $A$-facet $H$, and $\mathcal R=\mathcal N\setminus\{H\}$. Put $E=F_i$ if $\mathcal M\in\mathcal Q_i$, and $E=F$ if $\mathcal M\in\mathcal Q_{A'}$ with $A\subsetneq A'$. Then $E\cap H=A$ and $\mathcal R$ avoids $E\cup H$, so by (C1)$_p$ (or directly, if $|A|=d$) there is an irredundant proper chain $H,H_1,\ldots,E$ inside $E\cup H$, and by Lemma~\ref{lem:chain} $H_1=(H\setminus\{a\})\cup\{b\}$ with $b\in E\setminus H$ and $a\notin F$. Let $\mathcal N'=\mathcal R\cup\{H_1\}$. If $b$ is the free vertex of $E$, then $H_1=E$ and $\mathcal N'$ is in a clone block. Otherwise $b\in C\setminus A$, and $H_1$ is an $(A\cup\{b\})$-facet, so $\mathcal N'\in\mathcal Q_{A\cup\{b\}}$, an earlier block. Finally $b\in V_{\mathcal M}$: this is clear if $E=F_i$; if $E=F$, then $b\neq f$ (otherwise $H_1=F$ and $A=C$, while $A\subsetneq A'\subseteq C$), and since $A\cup\{b\}$ and $A'$ are comparable sets properly containing $A$, we get $b\in A'\subseteq V_{\mathcal M}$.

Let $\mathcal N\in\mathcal Q_0$. Then $\nu(\Delta)=p+1$ and $\mathcal L=\mathcal N\cup\{F\}$ is a $(p+1)$-matching. For $P\in\mathcal N$, condition (C2)$_p$ gives a path from $\mathcal N=\mathcal L\setminus\{F\}$ to $\mathcal L\setminus\{P\}$ in $G_p(\Delta;V_{\mathcal L})$; let $\mathcal N'$ follow $\mathcal N$ on it. Then $V_{\mathcal N'}=(V_{\mathcal N}\setminus\{a\})\cup\{b\}$ with $b\in F$, and the facet of $\mathcal N'$ meeting $F$ meets it in $\{b\}$, so $b\neq f$ and $\mathcal N'$ lies in a clone block or in $\mathcal Q_{\{b\}}$. If $\mathcal M\in\mathcal Q_i$, then $b\in C\subseteq V_{\mathcal M}$; if $\mathcal M\in\mathcal Q_{A'}$, then $\{b\}$ and $A'$ are comparable, so $b\in A'$.
\end{proof}


For the highest power, put $q=\nu(\Delta)$. Condition (C2)$_q$ is vacuous, and (C1)$_q$ can be stated without the restriction $|F\cap G|>0$: a $(q-1)$-matching avoiding two disjoint facets would give a $(q+1)$-matching.

\begin{Corollary}\label{cor:top}
Let $\Delta$ be a pure simplicial forest of dimension $d\geq 1$ and $q=\nu(\Delta)$. Then $I(\Delta)^{[q]}$ has linear quotients if and only if it is linearly related, if and only if $\nu(\Delta_{V\setminus(F\cup G)})<q-1$ for all facets $F,G$ with $|F\cap G|<d$ that are not connected in codimension one inside $\Delta_{F\cup G}$.
\end{Corollary}



\section{Squarefree powers of path ideals}\label{sec:path-ideals}
In this section, we study the linear quotients of squarefree powers of path ideals of trees.
For directed paths of rooted trees, we apply
Theorem~\ref{thm:main} to obtain a combinatorial
characterization. For undirected paths of forests, the
associated path complex need not be a simplicial forest,
so a different approach is required.


\subsection{Directed paths in rooted trees}\label{subsec:rooted-trees}

Let $\Gamma$ be a rooted tree, with edges directed away from the root,
and set $V=V(\Gamma)$.
For $t\geq2$, let $\Delta=\Gamma_t$ be the simplicial complex whose
facets are the vertex sets of the directed paths with $t$ vertices.
If nonempty, $\Delta$ is a pure simplicial tree of dimension
$t-1$~\cite[Theorem~2.7]{JV}.
For $p\geq1$, consider the following conditions:

\begin{enumerate}
\item[\textnormal{(R1)$_p$}] $\nu(\Delta_{V\setminus(F\cup G)})<p-1$ whenever $F,G$ are directed $t$-paths with the same initial vertex and $|F\cap G|<t-1$, as in Figure~\ref{fig:rooted-branch};
\item[\textnormal{(R2)$_p$}] the undirected subgraph of $\Gamma$ induced by $V_{\mathcal L}$ is connected for every $(p+1)$-matching $\mathcal L$ of $\Delta$.
\end{enumerate}

\begin{figure}[ht]
\centering
\begin{tikzpicture}[x=1cm,y=1cm,>=stealth,every node/.style={font=\small,inner sep=2pt}]
\node (a) at (0,0) {$x_1=y_1$};
\node (b) at (1.35,0) {$\cdots$};
\node (c) at (2.65,0) {$x_i=y_i$};
\node (d) at (4.25,.38) {$x_{i+1}$};
\node (e) at (5.35,.38) {$\cdots$};
\node (f) at (6.4,.38) {$x_t$};
\node (g) at (4.25,-.38) {$y_{i+1}$};
\node (h) at (5.35,-.38) {$\cdots$};
\node (j) at (6.4,-.38) {$y_t$};
\draw[->] (a)--(b); \draw[->] (b)--(c);
\draw[->] (c)--(d); \draw[->] (d)--(e); \draw[->] (e)--(f);
\draw[->] (c)--(g); \draw[->] (g)--(h); \draw[->] (h)--(j);
\end{tikzpicture}
\caption{A pair in \textnormal{(R1)$_p$}, with $1\leq i\leq t-2$.}
\label{fig:rooted-branch}
\end{figure}


\begin{Lemma}\label{lem:rooted-C1}
The conditions \textnormal{(C1)$_p$} and \textnormal{(R1)$_p$} are equivalent for $\Delta$.
\end{Lemma}

\begin{proof}
Assume (C1)$_p$. Paths with the same initial vertex cannot meet again after separating. Thus a pair $F,G$ in (R1)$_p$ has a common initial segment, and the only $t$-paths contained in $F\cup G$ are $F,G$; since $|F\cap G|<t-1$, they are not connected in codimension one inside $\Delta_{F\cup G}$. As $0<|F\cap G|$, condition (C1)$_p$ gives the required inequality.

Conversely, assume (R1)$_p$ and let $F,G$ be a pair occurring in (C1)$_p$. Their initial vertices are ancestors of any common vertex, so they are comparable. After interchanging $F,G$, the initial vertex of $G$ lies on $F$. Follow $F$ to that vertex and then follow $G$, obtaining a directed path $Q\subseteq F\cup G$. Its first $t$ vertices form a facet $H$, and its consecutive $t$-vertex subpaths give a proper chain from $H$ to $G$. Since $F,G$ are not connected inside $\Delta_{F\cup G}$, we have $H\neq F$ and $|F\cap H|<t-1$. The pair $F,H$ has the same initial vertex, so (R1)$_p$ gives $\nu(\Delta_{V\setminus(F\cup G)})\leq\nu(\Delta_{V\setminus(F\cup H)})<p-1$.
\end{proof}

\begin{Lemma}\label{lem:rooted-C2}
If $\Delta$ satisfies \textnormal{(R1)$_p$} and \textnormal{(R2)$_p$}, then it satisfies \textnormal{(C2)$_p$}.
\end{Lemma}

\begin{proof}
Fix a $(p+1)$-matching $\mathcal L$. Contracting its paths in the connected graph of (R2)$_p$ gives a connected graph on $\mathcal L$; hence it suffices to join the deletions of two paths $E,E'$ linked by an edge. Relabel them so that this edge is directed from $x_j\in E=\{x_1\to\cdots\to x_t\}$ to $w\in E'$. The vertex $w$ is the initial vertex of $E'$: otherwise its predecessor in $E'$ and $x_j$ would be distinct parents. Write $E'=\{w_1\to\cdots\to w_t\}$, with $w_1=w$, and put $\mathcal R=\mathcal L\setminus\{E,E'\}$.

Let $H$ consist of the first $t$ vertices of $Q=(x_1\to\cdots\to x_j\to w_1\to\cdots\to w_t)$. If $j\leq t-2$, then $E,H$ have the same initial vertex and $|E\cap H|=j$, while the $(p-1)$-matching $\mathcal R$ avoids $E\cup H$, contradicting (R1)$_p$. Thus $j\geq t-1$, and $H$ equals $E$ or shares $t-1$ vertices with it. The consecutive $t$-vertex subpaths of $Q$ join $H$ to $E'$ inside $E\cup E'$. Adjoining $\mathcal R$ to each facet of the resulting proper chain from $E$ to $E'$ gives a path from $\mathcal L\setminus\{E'\}$ to $\mathcal L\setminus\{E\}$ in $G_p(\Delta;V_{\mathcal L})$.
\end{proof}



\begin{Corollary}\label{cor:rooted-trees}
Let $\nu=\nu(\Gamma_t)$ and $1\leq p\leq\nu$. Then $I(\Gamma_t)^{[p]}$ has linear quotients if and only if it has a linear resolution, if and only if it is linearly related, if and only if \textnormal{(R1)$_p$} and \textnormal{(R2)$_p$} hold. In this case $\nu-1\leq\nu_0(\Gamma_t)\leq p$, so $p\in\{\nu-1,\nu\}$.
\end{Corollary}

\begin{proof}
The two lemmas and Theorem~\ref{thm:main} prove sufficiency and the numerical assertions. If $I(\Delta)^{[p]}$ is linearly related, that theorem and Lemma~\ref{lem:rooted-C1} give (R1)$_p$. For every $(p+1)$-matching $\mathcal L$, Proposition~\ref{prop:C2} gives codimension one connectedness of $\Delta_{V_{\mathcal L}}$. Its facets are connected subgraphs of $\Gamma$, consecutive facets of a proper chain intersect, and their union is $V_{\mathcal L}$ since they include $\mathcal L$; hence (R2)$_p$ holds.
\end{proof}

In particular, $I(\Gamma_t)^{[\nu]}$ has linear quotients exactly when
\textnormal{(R1)$_\nu$} holds, since \textnormal{(R2)$_\nu$} is vacuous.


\subsection{Undirected \texorpdfstring{$3$}{3}-paths in forests}\label{sec:3path}

We next consider path ideals generated by all undirected paths of a fixed length in a forest. The corresponding path complex need not be a simplicial forest. For example, the three $3$-paths in the complete bipartite graph $K_{1,3}$ form a
subcollection with no leaf. Therefore Theorem~\ref{thm:main} does not apply here. We prove directly that the highest nonzero squarefree power of the $3$-path ideal of forests admit linear quotients.

All graphs are finite and simple. A \emph{$t$-path} is a path on $t$ distinct vertices, $G^{(t)}$ is the hypergraph of the vertex sets of the $t$-paths of $G$, and $I_t(G)=I(G^{(t)})$. For $U\subseteq V(G)$, $G\setminus U$ is obtained by removing the vertices of $U$, and $N_G(u)$ is the set of neighbours of $u$. For a hypergraph $\mathcal H$ we denote by $\mathcal M_k(\mathcal H)$ its set of $k$-matchings, with $\mathcal M_0(\mathcal H)=\{\emptyset\}$, and as before $V_M=\bigcup_{F\in M}F$.

We use the following formulation of the linear quotient criterion
for sets of matchings.

\begin{Definition}\label{def:matching-order}
Let $\mathcal Q$ be a set of matchings of a hypergraph. A total order $>$ on $\mathcal Q$ is a \emph{linear quotient order} if for all $M>N$ in $\mathcal Q$ there is $W>N$ in $\mathcal Q$ with $|V_W\setminus V_N|=1$ and $V_W\setminus V_N\subseteq V_M\setminus V_N$. For disjoint sets of matchings, $\mathcal Q_1>\mathcal Q_2$ means that every element of $\mathcal Q_1$ comes before every element of $\mathcal Q_2$.
\end{Definition}

\begin{Remark}\label{rem:lqo-ideal}
If $\mathcal H$ is uniform and $\mathcal M_\nu(\mathcal H)$, $\nu=\nu(\mathcal H)\geq 1$, has a linear quotient order, then distinct maximum matchings have distinct vertex sets (for $M>N$ the set $V_M\setminus V_N$ is nonempty), so by the criterion stated at the beginning of Section~\ref{sec:LQ} the monomials $x_{V_M}$, in this order, give linear quotients of $I(\mathcal H)^{[\nu]}$. Moreover, if $\mathcal Q$ has a linear quotient order and $F$ is an edge disjoint from all the sets $V_Q$, $Q\in\mathcal Q$, then $\{Q\cup\{F\}:Q\in\mathcal Q\}$ has a linear quotient order, since $V_{Q\cup\{F\}}=V_Q\sqcup F$.
\end{Remark}

\begin{Lemma}\label{lemma.lq}
Let $G=G_1\sqcup\cdots\sqcup G_r$ be a disjoint union of graphs, and set
$\nu_i=\nu(G_i^{(3)})$ and $\nu=\nu(G^{(3)})$.
If every $\mathcal M_{\nu_i}(G_i^{(3)})$ has a linear quotient order,
then $\nu=\nu_1+\cdots+\nu_r$ and
$\mathcal M_\nu(G^{(3)})$ has a linear quotient order.
\end{Lemma}

\begin{proof}
A $3$-path lies in one component, so a maximum matching is $M=M_1\sqcup\cdots\sqcup M_r$ with $M_i$ a maximum matching of $G_i^{(3)}$. Order these matchings lexicographically: $M>N$ if $M_q>N_q$ for the first $q$ with $M_q\neq N_q$. If $M>N$ and $W_q$ is the matching given by the order on $\mathcal M_{\nu_q}(G_q^{(3)})$ for $M_q>N_q$, replacing $N_q$ by $W_q$ in $N$ gives the required $W$.
\end{proof}

Next, we prove the main result of this section.
\begin{Theorem}\label{thm:3path}
Let $G$ be a forest and $\nu=\nu(G^{(3)})$. Then $\mathcal M_\nu(G^{(3)})$ has a linear quotient order. Hence, if $\nu\geq 1$, then $I_3(G)^{[\nu]}$ has linear quotients and a $3\nu$-linear resolution.
\end{Theorem}

\begin{proof}
By Lemma~\ref{lemma.lq} we may assume that $G$ is a tree, and we use induction on $n=|V(G)|$; for $n\leq 3$ there can be only one maximum matching. By the induction hypothesis and Lemma~\ref{lemma.lq}, the statement holds for all forests with fewer than $n$ vertices. We shall often use the following observation: if a $3$-path $F$ lies in some maximum matching, then the maximum matchings containing $F$ are the sets $Q\cup\{F\}$ with $Q\in\mathcal M_{\nu-1}((G\setminus F)^{(3)})$, so they have a linear quotient order by induction and Remark~\ref{rem:lqo-ideal}.

Let $u_1,u_2,\ldots,u_\ell$ be a longest path of $G$ and $\mathcal J=\mathcal M_\nu(G^{(3)})$. Then $u_1$ is a leaf, and so is every neighbour of $u_2$ other than $u_3$.

\emph{Case I: some $\widetilde M\in\mathcal J$ avoids $u_1$ and $u_2$.} Then $N_G(u_2)=\{u_1,u_3\}$, since another neighbour $z$ would be a leaf not covered by $\widetilde M$, and we could add $\{u_1,u_2,z\}$ to $\widetilde M$. Let $F=\{u_1,u_2,u_3\}$, the only $3$-path through $u_1$, and split $\mathcal J$ into
\[
\mathcal J_1=\{M:u_1,u_2\notin V_M\},\quad\mathcal J_2=\{M:u_2\in V_M,\ u_1\notin V_M\},\quad\mathcal J_3=\{M:F\in M\}.
\]
The set $\mathcal J_1$ is $\mathcal M_\nu((G\setminus\{u_1,u_2\})^{(3)})$, which has a linear quotient order by induction, and $\mathcal J_3$ has one by induction and Remark \ref{rem:lqo-ideal}. A matching in $\mathcal J_2$ contains $F_y=\{u_2,u_3,y\}$ for some $y\in N_G(u_3)\setminus\{u_2\}$, and the set $\mathcal J_2(y)$ of these matchings has a linear quotient order by induction and Remark \ref{rem:lqo-ideal} (note that $u_1$ is isolated in $G\setminus F_y$). We put the sets $\mathcal J_2(y)$ with $y\neq u_4$ first, in any order, and $\mathcal J_2(u_4)$ last.

Let $M\in\mathcal J_2(y)$ come before $N\in\mathcal J_2(y')$ with $y'\neq y$. Then $y\neq u_4$ is not on the path. Let $T$ be the component of $G\setminus u_3$ containing $y$. By maximality of the path, the vertices of $T$ other than $y$ are leaves adjacent to $y$, and there is at most one of them: if $a,b$ were two, then $(M\setminus\{F_y\})\cup\{F,\{a,y,b\}\}$ would be a $(\nu+1)$-matching. So every $3$-path through $y$ contains $u_3$. Since $u_3$ is covered by $F_{y'}\in N$, we get $y\notin V_N$, and $W=(N\setminus\{F_{y'}\})\cup\{F_y\}$ satisfies $V_W\setminus V_N=\{y\}\subseteq V_M\setminus V_N$. Hence $\mathcal J_2$ has a linear quotient order.

We order $\mathcal J_3>\mathcal J_2>\mathcal J_1$. If $M\in\mathcal J_3$ and $N\in\mathcal J_2$, replacing $F_y\in N$ by $F$ gives $W\in\mathcal J_3$ with $V_W\setminus V_N=\{u_1\}$. If $M\in\mathcal J_2\cup\mathcal J_3$ and $N\in\mathcal J_1$, then $u_3\in V_N$ by maximality; if $E\in N$ contains $u_3$ and $x\in E\cap N_G(u_3)$, replacing $E$ by $\{u_2,u_3,x\}$ gives $W\in\mathcal J_2$ with $V_W\setminus V_N=\{u_2\}$.

\emph{Case II: every $M\in\mathcal J$ covers $u_2$, and some $\widetilde M\in\mathcal J$ avoids $u_1$.} Let $\mathcal J_1=\{M:u_1\notin V_M\}=\mathcal M_\nu((G\setminus u_1)^{(3)})$, which has a linear quotient order by induction, and $\mathcal J_2=\{M:u_1\in V_M\}$. A matching in $\mathcal J_2$ contains $F_x=\{u_1,u_2,x\}$ for some $x\in N_G(u_2)\setminus\{u_1\}$, and the set $\mathcal J_2(x)$ of these matchings has a linear quotient order by induction and Remark \ref{rem:lqo-ideal}. Put the sets $\mathcal J_2(x)$ with $x\neq u_3$ first and $\mathcal J_2(u_3)$ last. If $M\in\mathcal J_2(x)$ comes before $N\in\mathcal J_2(x')$, then $x$ is a leaf, not covered by $N$ because $u_2\in F_{x'}$, and $W=(N\setminus\{F_{x'}\})\cup\{F_x\}$ gives $V_W\setminus V_N=\{x\}$. Finally let $\mathcal J_2>\mathcal J_1$: for $N\in\mathcal J_1$, if $E\in N$ contains $u_2$ and $x\in E\cap N_G(u_2)$, then $W=(N\setminus\{E\})\cup\{\{u_1,u_2,x\}\}\in\mathcal J_2$ and $V_W\setminus V_N=\{u_1\}$.

\emph{Case III: every $M\in\mathcal J$ covers $u_1$.} Then $M$ contains some $\{u_1,u_2,x\}$, and $x$ does not depend on $M$: if $M$ contains $\{u_1,u_2,x\}$ and $N$ contains $\{u_1,u_2,y\}$ with $x\neq y$, we may assume $y\neq u_3$, so that $y$ is a leaf not covered by $M$, and replacing $\{u_1,u_2,x\}$ by $\{x,u_2,y\}$ in $M$ gives a maximum matching avoiding $u_1$. So all maximum matchings contain the same $3$-path and the conclusion follows by induction and Remark \ref{rem:lqo-ideal}.

Since a matching covering $u_1$ covers $u_2$, the three cases are exhaustive. The last assertion follows from Remark~\ref{rem:lqo-ideal}.
\end{proof}



The following examples show that Theorem~\ref{thm:3path} does not extend to $4$-path ideals of trees or to $3$-path ideals of arbitrary graphs.


\begin{Example}\label{ex:four-paths}
Let $G$ be the tree on $[7]$ with edges $\{1,2\}$, $\{2,3\}$, $\{3,4\}$, $\{4,5\}$, $\{3,6\}$, $\{6,7\}$. Every $4$-path contains $3$, so $\nu(G^{(4)})=1$ and
\[
I_4(G)^{[1]}=I_4(G)=x_3(x_1x_2x_4,\,x_1x_2x_6,\,x_2x_4x_5,\,x_4x_5x_6,\,x_2x_6x_7,\,x_4x_6x_7).
\]
The nonzero graded Betti numbers of $S/I_4(G)$ are $\beta_{0,0}=1$, $\beta_{1,4}=6$, $\beta_{2,5}=6$ and $\beta_{3,7}=1$. Thus $I_4(G)$ is linearly related, but
$\reg I_4(G)=5$, so it does not have a linear resolution.
\end{Example}


\begin{Example} \label{ex.arbG}
    Let $G$ be the fully whiskered cycle with  $V(G)=\{x_i,y_i\mid 1\leq i\leq 10\}$, and $$E(G)=\bigl\{\{x_i,x_{i+1}\}:1\leq i\leq 10\bigr\} \cup \bigl\{\{x_i,y_i\}:1\leq i\leq 10\bigr\},$$ where $x_{11}=x_1$. Since every $3$-path in $G$ contains at least two cycle vertices, we have $\nu(G^{(3)})\leq 5$. The five pairwise vertex-disjoint $3$-paths $\{y_{2i-1},x_{2i-1},x_{2i}\}$, $1\leq i\leq 5$, show that $\nu(G^{(3)})=5$. Since $\beta_{4,20}(I_3(G)^{[5]})=1$, we have $\reg(I_3(G)^{[5]})\geq 16>15$. Thus, $I_3(G)^{[5]}$ does not have a $15$-linear resolution.
\end{Example}

\end{document}